\documentclass{article}
\usepackage[utf8]{inputenc}

\title{Hamiltonicity of 3-tough $(K_2 \cup 3K_1)$-free graphs}
\author{ 
    {\large Elizabeth Grimm, Andrew Hatfield}\smallskip \vspace{-.2cm} \\ 
	\medskip  \normalsize{Department of Mathematics, Illinois State University, Normal, IL 61790}\vspace{-.33cm}\\ 
	{\normalsize {\tt evgrimm@ilstu.edu, abhatfi@ilstu.edu}}
}
\date{\today}
\usepackage[shortlabels]{enumitem}
\usepackage{amssymb}
\usepackage{amsfonts}
\usepackage{amsthm}
\usepackage{amsmath}
\usepackage{indentfirst}
\usepackage{mathtools}

\DeclarePairedDelimiter\floor{\lfloor}{\rfloor}
\newtheorem{THM}{Theorem}

\newtheorem{LEM}[THM]{Lemma}
\newtheorem*{theorem*}{Theorem}
\newtheorem{CLAIM}{Claim}
\newcommand{\pf}{\textbf{Proof}.\quad}
\usepackage[
pdfauthor={ESYZ},
pdftitle={},
pdfstartview=XYZ,
bookmarks=true,
colorlinks=true,
linkcolor=blue,
urlcolor=blue,
citecolor=blue,
bookmarks=false,
linktocpage=true,
hyperindex=true
]{hyperref}
\begin{document}
\newcommand{\oC}{\overset{\rightharpoonup }{C}}
\newcommand{\iC}{\overset{\leftharpoonup }{C}}
\providecommand{\keywords}[1]
{
  \small	
  \textbf{\textit{Keywords: }} #1
}

\maketitle


\begin{abstract}
Chv\'{a}tal conjectured in 1973 the existence
of some constant $t$ such that all $t$-tough graphs with at least three vertices are hamiltonian. While the conjecture
has been proven for some special classes of graphs, it remains open in general. We say that a graph
is $(K_2 \cup 3K_1)$-free if it contains no induced subgraph isomorphic to
$K_2 \cup 3K_1$, where $K_2 \cup 3K_1$ is the disjoint union of an edge and three isolated vertices. In this paper, we show that every 3-tough $(K_2 \cup 3K_1)$-free graph with at least three vertices is hamiltonian.\\
\end{abstract}


\section{Introduction}
Let $G$ be a simple graph and let $E(G), V(G)$ denote its edge and vertex set respectively. For 
$v \in V(G)$, denote by $N_G(v)$ the set of neighbors of $v$ in $G$, and let $d_G(v) = |N_G(v)|$ be the degree
of $v$ in $G$. For two disjoint subgraphs $H_1, H_2$ of $G$, $N_{H_1}(H_2)$ denotes the set of neighbors of vertices of $H_2$ in $G$ 
that are contained in $V(H_1)$. 
For any subset $S \subseteq V(G)$,  $G[S]$ is the subgraph of $G$ induced on $S$, 
$G - S$ denotes the subgraph $G[V(G) \setminus S]$, and $N(S) = \cup_{v \in S} N_G(v)$.
For any disjoint $A, B \subseteq V(G)$,  $E_G(A, B)$ is  the set of all edges  with one end-vertex
in $A$ and the other end-vertex in $B$.  If $u$ and $v$ are adjacent in $G$, we write $u\sim v$.

%
%

We say that $G$ is hamiltonian if there exists a cycle which contains every vertex of $G$.
We say that $G$ is $H$-free if there does not exist an induced copy of $H$ in $G$. Denote by $c(G)$ the number of components of $G$.
Let $t \geq 0$ be a real number. We say a graph $G$ is $t$-tough if for each cut set $S$ of $G$ we have $t \cdot c(G-S) \leq |S|$. The toughness of a graph $G$, denoted $\tau(G)$, is the maximum value of $t$ for which $G$ is $t$-tough if $G$ is non-complete, and is defined to be $\infty$ if $G$ is complete. Chv\'{a}tal introduced the notion of toughness in his 1973 paper \cite{MR316301}, where he also conjectured the existence of a constant $t$ such that every $t$-tough graph on at least three vertices is hamiltonian. The conjecture has been verified for certain special classes of graphs, but remains open in general. Recent work has proven the conjecture for $2K_2$-free graphs \cite{broersma2k2, shan2k2, ota2021toughness}, $(P_2 \cup P_3)$-free graphs \cite{shanp2up3}, $(K_2 \cup 2K_1)$-free graphs \cite{1tough}, and planar chordal graphs. We refer the reader to \cite{bauersurvey2} for a survey on more related results.


In this paper, we support Chv\'{a}tal's conjecture by proving the following result:

\begin{THM}\label{theorem1}
If $G$ is a 3-tough $(K_2 \cup 3K_1)$-free graph on at least 3 vertices, then $G$ is hamiltonian.
\end{THM}

The remainder of this paper is organized as follows: in Section \ref{preliminaries}, we discuss results necessary for the proof of Theorem \ref{theorem1} and in Section \ref{section3} we prove Theorem \ref{theorem1}.

\section{Preliminaries}\label{preliminaries}

In this section, we give results necessary to complete the proof of Theorem \ref{theorem1}. 


\begin{LEM}[Dirac \cite{hamconn1}, Ore \cite{hamconn2}]\label{1.5}
    Let $G$ be a graph on $n$ vertices such that $\delta(G) \geq \frac{n+1}{2}$. Then $G$ is hamiltonian-connected. 
\end{LEM}

\begin{LEM}[Bauer et al. \cite{MR1336668}]\label{bauer} 
 Let $G$ be a $t$-tough graph on $n \geq 3$ vertices with $\delta(G) > n/(t + 1) - 1$. Then $G$ is hamiltonian.
\end{LEM}

\begin{LEM}[Li et al. \cite{1tough}]\label{k2cup2k1}
    Let $R$ be an induced subgraph of $P_4$, $K_1 \cup P_3$ or $K_2 \cup 2K_1$. Then
every $R$-free 1-tough graph on at least three vertices is hamiltonian.
\end{LEM}

The following lemma is a consequence of Menger's theorem, which can be found in \cite{bondy}. For a positive integer $k$, define 
$[1,k] = \{1,2,\cdots, k\}$.


\begin{LEM}\label{3.2.6}
Let $G$ be a $k$-connected graph and $X_1,X_2$ be distinct subsets of $V(G)$. 
Then there exist $k$ internally disjoint paths $P_1, \ldots, P_k$ such that 
\begin{enumerate}[(a)]
    \item $|V(P_i) \cap X_1| = |V(P_i) \cap X_2| = 1$, and $P_i$ is internally disjoint from each $X_1$ and $X_2$.
    \item if $|X_i| \geq k$ for some $i \in [1,2]$, then $V(P_j) \cap X_i \ne V(P_\ell) \cap X_i$ for all distinct $j, \ell \in [1,k]$.
    \item if $|X_i| < k$ for some $i \in [1,2]$, then every vertex of $X_i$ is an end-vertex of some path $P_j$ for $j \in [1,k]$.
\end{enumerate}
\end{LEM}
The following lemma provides some structural properties of $(K_2 \cup 3K_1)$-free graphs. 

\begin{LEM}\label{lemma:cut-structure}
Let $G$ be a connected ($K_2 \cup 3K_1$)-free graph, and $S 
\subseteq V(G)$ be a cut set such that $G - S$ has at least three components. Then we have
the following statements:
\begin{enumerate}[(a)]
	\item If $G - S$ has a nontrivial component, then $G - S$ has exactly three components. \label{parta}
	\item If $G - S$ has a nontrivial component, then the component is ($K_2 \cup K_1$)-free.
\end{enumerate}
\end{LEM}

\pf   For  part $(a)$, 
let $D$ denote a nontrivial component of $G - S$. Assume for the sake of contradiction that $G - S$ has more than three  components.
Taking an edge from $D$ and a single vertex from three other components, respectively, gives an induced copy of $K_2 \cup 3K_1$.
This gives a contradiction to the ($K_2 \cup 3K_1$)-freeness of $G$. 

For part $(b)$,      $G - S$ must have exactly three  components by part (a). Assume for the sake of contradiction
    that the nontrivial component  is not $(K_2 \cup K_1)$-free. 
    Then taking an induced copy of $K_2 \cup K_1$ from this component and one vertex each from the other two components gives 
     an induced copy of $K_2 \cup 3K_1$, which contradicts 
the ($K_2 \cup 3K_1$)-freeness of $G$. 
\qed
\\

Note that in any $(K_2 \cup 3K_1)$-free graph $G$, the components yielded by any cut set $S$ such that $c(G-S) \geq 3$ must be $(K_2 \cup K_1)$-free. The following lemmas deal with the structure of $(K_2 \cup K_1)$-free graphs. Lemma $\ref{lemma:trivialcomponents}$ is 
used in the proof of Lemma \ref{mindegreeD}.

\begin{LEM}\label{lemma:trivialcomponents}
If $G$ is a ($K_2 \cup K_1$)-free graph and $S$ is a cut set of $G$, then 
every component of $G - S$ is trivial.
\end{LEM}

\pf Assume there exists some nontrivial component of $G - S$. Since $S$ is a cut set, 
it must disconnect $G$ into at least two components. Taking an edge from a nontrivial component
and a vertex from another component gives an induced copy of $K_2 \cup K_1$, contradicting the 
$(K_2 \cup K_1)$-freeness of $G$. 
\qed 

The independence number of a graph $G$, denoted $\alpha (G)$, is the size of a largest independent set of $G$.


\begin{LEM}\label{mindegreeD}
    Let $t > 0$ be real and $G$ be a $(K_2 \cup K_1)$-free graph on $n$ vertices with $\alpha(G) \leq \frac{n}{t+1}$. Then $\delta(G) \geq n - \frac{n}{t+1}$.
\end{LEM}

\pf
Assume  $\delta(G) < n-\frac{n}{t+1}$. 
Let $v\in V(G)$ with $d_G(v)=\delta(G)$, and let 
$W=V(G)\setminus N_G(v)$.  Then $|W| =|V(G)|-\delta(G) > \frac{n}{t+1}$.
As $G$ is $(K_2\cup K_1)$-free, and $N_G(v)$
is a cut set of $G$, every component of $G-N_G(v)$ is trivial by Lemma~\ref{lemma:trivialcomponents}. Since $W=G-N_G(v)$,  $W$ is an independent set of $G$. However, $|W| > \frac{n}{t+1} \geq \alpha(G)$, giving a contradiction.
\qed


\section{Proof of Theorem \ref{theorem1}}\label{section3}
For any $u,v \in V(G)$, we call a path $P$ connecting $u$ and $v$ a $(u,v)$-path. 
Let $P$ be an $(x,y)$-path and $Q$ be a $(y,z)$-path such that $V(P) \cap V(Q) = \{y\}$. Then $xPyQz$ denotes a path from $x$ to $z$. 
If $P$ is an $(x,y)$-path and $Q$ is a disjoint $(z, w)$-path where $y \sim z$, we denote by $xPyzQw$ the concatenation 
of $P$ and $Q$ through the edge $yz$.  
Let $\oC$ 
be a cycle with some fixed orientation. For any $u,v \in V(C)$, we denote by $u\oC v$ the path from $u$ to $v$ following the orientation of $C$. Similarly, we denote by $u \iC v$ the inverse path from $u$ to $v$. 
The immediate successor of $u$ on $\oC$ is denoted by $u^+$. 

\proof[Proof of Theorem~\ref{theorem1}]
Let $G$ be a 3-tough $(K_2 \cup 3K_1)$-free graph. We may assume that $G$ is not complete, otherwise there exists a hamiltonian cycle. Therefore, $G$ is 6-connected. By Theorem \ref{bauer},
we may assume $\delta(G) \leq \frac{n}{4} -1$. Since $\delta(G) \geq 6$, we get $n \geq 28$.
Let $C$ be a longest cycle of $G$.

\begin{CLAIM}\label{claim1}
$|V(C)| \geq \frac{3n}{4}$.
\end{CLAIM}

\pf We assume first that there exist $u,v \in V(G)$ with $u \not \sim v$ such that $|N(u) \cup N(v)| \leq \frac{n}{4}$. Let $S = N(u) \cup N(v)$.

Note that the components of $G - S$ cannot all be trivial, as this would imply $\frac{|S|}{c(G-S)} \leq \frac{\frac{n}{4}}{\frac{3n}{4}} < 3$, but 
$G$ is 3-tough. Thus, there must exist a nontrivial component $D$ in $G - S$. Since each of $u,v$ are components of $G - S$ and $G - S$ has a nontrivial component, it follows that $c(G - S) \geq 3$. Thus $G - S$ has exactly 3 components by Lemma \ref{lemma:cut-structure}.
Since $G$ is 3-tough, it follows that $\alpha(G) \leq \frac{n}{4}$.

Also, we know that $\delta(D) \geq |V(D)| - \frac{n}{4} \geq \frac{|V(D)| + 1}{2}$, as $n \geq 12$ and $|V(D)| \geq \frac{3n}{4}-2$ by Lemma \ref{mindegreeD}. Thus $D$ is hamiltonian-connected by Lemma \ref{1.5}. Since $G$ is 2-connected, by Lemma \ref{3.2.6} we can find in $G$ two disjoint paths $P_1$ from $u$ to some $x_1 \in V(D)$, and $P_2$ from $u$ to some $x_2 \in V(D)$ where $x_1 \ne x_2$, and each $P_i$ is internally disjoint from $D$. As $D$ is hamiltonian-connected, we can find a hamiltonian path $Q$ in $D$ from $x_1$ to $x_2$.
Define $C' = uP_1x_1Qx_2P_2u$. Then $|V(C')| \geq |V(D)| + 3 \geq \frac{3n}{4} - 2 + 3 = \frac{3n}{4}+1 > \frac{3n}{4}$. 
As $C$ is a longest cycle of $G$, $|V(C)| \geq |V(C')| \geq \frac{3n}{4}$.

We then assume that for any $u,v \in V(G)$ with $u \not\sim v$, it holds that 
$|N(u) \cup N(v)| > \frac{n}{4}$. Note that $\delta(G) \leq \frac{n}{4} - 1$ by our earlier assumption. Let $u \in V(G)$ with $d_G(u) = \delta(G)$.
Define $G_1 = G - (N_G(u) \cup \{u\})$.
We know that $G_1$ is $(K_2 \cup 2K_1)$-free, as there are no edges between $G_1$ and $u$ and the original graph is $(K_2 \cup 3K_1)$-free.
If $G_1$ is 1-tough, then it has a hamiltonian cycle by Lemma \ref{k2cup2k1} which has at least $\frac{3n}{4}$ vertices, thus 
$|V(C)| \geq \frac{3n}{4}$. Thus we may assume $G_1$ is not 1-tough. Let $W$ be a tough set of $G_1$, i.e. $W$ is a cut set of $G_1$ such that $\frac{|W|}{c(G_1-W)} = \tau(G_1)$.

We claim below that $c(G_1-W) = 2$.
We first note that $G_1 - W$ has at least one nontrivial component. Otherwise, all components of $G_1 - W$ are trivial,
i.e. $c(G_1 - W) = |V(G_1)| - |W|$. Since $\tau(G_1) < 1$, we have $|W| < c(G_1 - W)$. Then $c(G_1 - W) > \frac{1}{2}|V(G_1)| \geq 
\frac{1}{2}(n - \frac{n}{4}) = \frac{3n}{8}$. Let $S = N(u) \cup \{u\} \cup W$. Then $c(G - S) = c(G_1-W)$ and $\frac{|S|}{c(G-S)} = \frac{|S|}{c(G_1-W)} \leq \frac{5n/8}{3n/8} < 3$, which is a contradiction as $G$ is 3-tough. Thus $G_1 - W$ has a nontrivial component. 

Next, assume that $G_1 - W$ has more than 2 components. Then by the argument above, at least one of the components is nontrivial. Taking an edge from a nontrivial component
and a vertex from each of the two others gives a copy of $K_2 \cup 2K_1$, which contradicts the $(K_2 \cup 2K_1)$-freeness of $G_1$. Thus, we must have 
$c(G_1 - W) = 2$. 

Let $D_1, D_2$ be the two components of $G_1 - W$.
Then as $G_1$ is not 1-tough, we get $|W| = 1$. Note that $|V(D_i)| \geq 2$. Otherwise, say $|V(D_2)| = 1$ and let $V(D_2) = \{v\}$, then $|N(u) \cup N(v)| \leq \frac{n}{4} - 1 + 1 = \frac{n}{4}$ which contradicts the assumption that 
$|N(u) \cup N(v)| > \frac{n}{4}$ for any nonadjacent $u,v \in V(G)$.
Since both $D_1$ and $D_2$ are nontrivial and $G_1$ is $(K_2 \cup 2K_1)$-free, each $D_i$ is a complete graph.

Then by Lemma \ref{3.2.6}, we can find in $G$ two disjoint paths $P_1$ from some $x_1 \in V(D_1)$ to some $x_2 \in V(D_2)$, and $P_2$ from some $y_1 \in V(D_1)$ to some $y_2 \in V(D_2)$ where $x_1 \ne y_1, x_2 \ne y_2$. As each $D_i$ is complete, we can find a hamiltonian path $Q_i$ in $D_i$ from $x_i$ to $y_i$. Then the cycle 
$$C' = x_1P_1x_2Q_2y_2P_2y_1Q_1x_1$$
satisfies $|V(C')| \geq |V(D_1)| + |V(D_2)| + 2 \geq \frac{3n}{4}$. As $C$ is a longest cycle, we have $|V(C)| \geq |V(C')| \geq \frac{3n}{4}$. $\hfill\blacksquare$
\\

Assume that $C$ is not hamiltonian, as otherwise we are done. Thus $G - V(C)$ has components. Orient $C$ in the clockwise direction and denote the orientation by $\oC$.

\newcommand{\xip}{\{x_i^+ \ | \ 1 \leq i \leq 6\}}

\begin{CLAIM}\label{hproperties}
Let $H$ be any component of $G - V(C)$. Then we have the following statements:
\begin{enumerate}[(a)]
    \item $|N_C(H)| \geq 2\tau(G) \geq 6$.
    \item for any two $x, y \in N_C(H)$, $xy \not \in E(C)$.
    \item for any two $x, y \in N_C(H)$, $x^+y^+ \not \in E(G)$. 
    \item $H$ is a trivial component.
\end{enumerate}
\end{CLAIM}

\pf
Let $H$ be a component of $G - V(C)$, and $x,y \in N_C(H)$.
Note that since $G$ is 3-tough, $2\tau(G) \geq 6$. 

For part $(a)$, assume $|N_C(H)| < 2\tau(G)$.  Then $\frac{|N_C(H)|}{c(G - N_C(H))} < \frac{2\tau(G)}{2} = \tau(G)$, contradicting the toughness of $G$. Thus we have $|N_C(H)| \geq 2\tau(G) \geq 6$.

For part $(b)$, $|N_C(H)| \geq 6$ by part $(a)$. If there exist distinct $x,y \in N_C(H)$ such that $xy \in E(C)$, then let $h_1 \in N_H(x), h_2 \in N_H(y)$.
Assume without loss of generality that $y = x^+$.
As $H$ is connected, there exists some $(h_1, h_2)$-path $P$ in $H$. Then the cycle $C' = xh_1Ph_2y\oC x$ is a cycle longer than $C$,
contradicting the maximality of $C$.

For part $(c)$, assume for the sake of contradiction that $x^+y^+ \in E(G)$. Let $h_1 \in N_H(x), h_2 \in N_H(y)$. 
Assume without loss of generality that $x$ appears before $y$ on the cycle. 
Again, as $H$ is connected, there exists some $(h_1, h_2)$-path $P$ in $H$. Then the cycle $C' = x h_1 P h_2 y \iC x^+ y^+ \oC x$ is a cycle
longer than $C$, contradicting the maximality of $C$. 

For part $(d)$, note that $\{x^+ \ | \ x \in N_C(H)\}$ is an independent set of $G$ by Claim \ref{hproperties} $(c)$. We assume that $H$ is nontrivial, then taking an edge
from $H$ and three vertices from the independent set $\{x^+ \ | \ x \in N_C(H)\}$ gives an induced copy of $(K_2 \cup 3K_1)$.
$\hfill\blacksquare$




\begin{CLAIM}\label{components}
$c(G - V(C)) \leq 3$. 
\end{CLAIM}

\pf
By Claim $\ref{hproperties}$ $(d)$, each component of $G - V(C)$ is trivial. For the sake of contradiction, assume $G - V(C)$ has at least 4 components, and let $x,y,z,w$ be 4 of them. 
Then the set $S = V(C) \setminus N_C(\{x,y,z,w\})$ is independent. Otherwise, taking an edge
from the set $S$ and three of $\{x,y,z,w\}$ gives an induced copy of $K_2 \cup 3K_1$. 
Then there must exist
3 consecutive vertices $u_1, u_2, u_3 \in V(C)$. Otherwise $|S| \geq \frac{n}{4}$, therefore the independent set $S \cup \{x,y,z,w\}$
has size at least $\frac{3n/4 + 4}{3} > \frac{n}{4}$, contradicting $\tau(G) \geq 3$. Note that none of $\{x,y,z,w\}$ 
can be adjacent to two consecutive $u_i$, as otherwise we may easily extend the cycle $C$. 
Without loss of generality, assume that $u_2 = u_1^+, x \sim u_1, y \sim u_2, z \sim u_3$. 
Then we must have $w \sim u_1$, as otherwise taking the edge $xu_1$ and the three vertices
$\{y,z,w\}$ gives an induced copy of $K_2 \cup 3K_1$. Similarly, we must have 
$w \sim u_2$. Thus the cycle $C' = u_1 w u_2 \oC u_1$ is longer than $C$, which is a contradiction. 
Therefore $c(G - V(C)) \leq 3$. $\hfill\blacksquare$

\begin{CLAIM}\label{claimwz}
Let $H$ be any component of $G - V(C)$, and $w,z \in N_C(H)$ be any two distinct vertices. Then for any $w_1 \in N_C(w^+)$, 
\begin{enumerate}[(a)]
    \item if $w, w_1, z$ appears in the order $w, z, w_1$ along $\oC$, then
    $z^+ \not\sim w_1^+$.
    \item if $w, w_1, z$ appear in the order $w, w_1, z$ along $\oC$, then
    $z^+ \not\sim w_1^-$.
\end{enumerate}
\end{CLAIM}

\pf Let $h_1 \in N_H(w), h_2 \in N_H(z)$. As $H$ is connected, we can find in $H$ an $(h_2, h_1)$-path $P$. For part (a), 
assume $z^+ \sim w_1^+$. Then define $$C' = w \iC w_1^+ z^+ \oC w_1 w^+ \oC z h_2 P h_1 w.$$ For part (b), assume $z^+ \sim w_1^-$. Then define $$C' = w \iC z^+ w_1^- \iC w^+ w_1 \oC z h_2 P h_1 w.$$ In either case, $C'$ is a longer cycle than $C$, contradicting the maximality of $C$. $\hfill\blacksquare$

By Claim $\ref{hproperties}$ $(d)$, we let $x \in V(G) \setminus V(C)$ be a component of $G - V(C)$. By Claim 
$\ref{hproperties}$ $(a)$, $|N_C(x)| \geq 2\tau(G) \geq 6$. We let $x_1, x_2, \cdots, x_6 \in N_C(x)$ be all distinct vertices
and assume they appear in the order $x_1, \cdots, x_6$ along $\oC$. Note that $\{x_1^+, \cdots, x_4^+\}$ is an independent set in $G$ by Claim \ref{hproperties} $(c)$.
Assume without loss of generality that $|V(x_4 \oC x_1)| \geq (|V(C)|-2)/2$. Label the 
vertices on this segment by $x_4, x_4^+, y_1, y_2, \cdots, y_t, x_1$. Since $x_4^+ \sim y_1$, 
we know that $|N(y_1) \cap \{x_1^+, \cdots x_4^+\}| \geq 2$, as otherwise taking the edge 
$x_4^+y_1$ and three vertices from the independent set $\{x_1^+, \cdots, x_4^+\}$ gives an induced copy of 
$K_2 \cup 3K_1$. We claim that $|N(y_2) \cap \{x_1^+, \cdots, x_4^+\}| = 0$. Otherwise,
by the same argument as above, we have $|N(y_2) \cap \{x_1^+, \cdots, x_4^+\}| \geq 2$.

If $x_1^+ \in (N(y_1) \cap N(y_2) \cap \{x_1^+, \cdots, x_4^+\})$, then there exists $x_j^+$
with $2 \leq j \leq 4$ such that $x_j^+ \sim y_2$, contradicting Claim $\ref{claimwz}$. 

If $x_4^+ \in N(y_1) \cap N(y_2) \cap \{x_1^+, \cdots, x_4^+\})$, then similarly there exists $x_j^+$
with $1 \leq j \leq 3$ such that $x_j^+ \sim y_1$, contradicting Claim $\ref{claimwz}$. 

If $(N(y_1) \cap N(y_2) \cap \{x_1^+, \cdots, x_4^+\}) \cap \{x_2^+, x_3^+\} \ne \emptyset$, we may assume that $N(y_1) \cap \{x_1^+, \cdots, x_4^+\} = \{x_4^+, x_j^+\}$
and $N(y_2) \cap \{x_1^+, \cdots, x_4^+\} = \{x_1^+, x_j^+\}$, otherwise we can apply Claim \ref{claimwz}. Then if $x \sim y_1$, 
the cycle $C' = x_j x y_1 \iC x_j^+ y_2 \oC x_j$is longer than $C$. Thus $x \not \sim y_1$. Then taking the independent set 
$\{x, x_1^+, x_2^+, x_3^+\} \setminus \{x_j^+\}$ and the edge $y_1x_j^+$ gives an induced copy of $K_2 \cup 3K_1$, which is a
contradiction. 

Thus $N(y_1) \cap N(y_2) \cap \{x_1^+, \cdots, x_4^+\} = \emptyset$. We may assume that 
$N(y_1) \cap \{x_1^+, \cdots, x_4^+\} = \{x_3^+, x_4^+\}$ and $N(y_2) \cap 
\{x_1^+, \cdots, x_4^+\} = \{x_1^+, x_2^+\}$, else we can apply Claim \ref{claimwz} as in previous arguments.
Note that if $x \sim y_1$, the cycle $C' = x_1 x y_1 \iC x_1^+ y_2 \oC x_1$ is longer than $C$. Thus 
$x \not \sim y_1$. 
Then taking the set $\{x, x_1^+, x_2^+\}$ and the edge $y_1x_3^+$ gives an induced copy
of $K_2 \cup 3K_1$, which is a contradiction. Thus $N(y_2) \cap \{x_1^+, x_2^+, \cdots, x_4^+\} = \emptyset$.

We claim that $|N(y_3) \cap \{x_1^+, \cdots, x_4^+\}| \geq 2$. Otherwise, taking 3 vertices from the independent set
$\{x_1^+, \cdots, x_4^+\}$ and the edge
$y_2y_3$ gives an induced copy of $K_2 \cup 3K_1$. Recall that the last vertex before $x_1$ on the segment
$x_4 \oC x_1$ is labelled by $y_t$. Define $Y$ as the set of all even-indexed $y_i$. Then $Y = \{y_2, y_4, \cdots, y_t\}$ in the 
case that $t$ is even and $Y = \{y_2, y_4, \cdots, y_{t-1}\}$ in the case that $t$ is odd. By applying a similar
argument as in the $y_2$ case to all even-indexed $y_i$, we see that $N(Y) \cap \{x_1^+, \cdots, x_4^+\} = \emptyset$. 

Note that $Y$ is an independent set of $G$, as otherwise 
taking an edge with end-vertices in $Y$ and 3 vertices from the set $\{x_1^+, \cdots, x_4^+\}$ gives
an induced copy of $K_2 \cup 3K_1$. Thus, $Y \cup \{x_1^+, \cdots, x_4^+\}$ is also an independent set in $G$.


By Claim $\ref{hproperties}$ $(d)$ and Claim $\ref{components}$, we have $|V(C)| \geq n - 3$. 
Thus $|Y| \geq \floor{\frac{1}{2}(\frac{|V(C)|-2}{2})} \geq \floor{\frac{n-5}{4}}$. 
Therefore $|Y \cup \{x_1^+, \cdots, x_4^+\}| > \frac{n}{4}$, contradicting $\tau(G) \geq 3$. This completes the proof.
\qed

\section*{Acknowledgements}
The authors wish to thank Dr. Songling Shan for her guidance during the research process and her valuable
feedback.

\bibliographystyle{plain}
\bibliography{bib}

\end{document}